\documentclass{article}
\usepackage{amssymb}
\usepackage{amsfonts}
\usepackage{amsmath}

\newtheorem{theorem}{Theorem}
\newtheorem{acknowledgement}[theorem]{Acknowledgement}

\newtheorem{definition}[theorem]{Definition}
\newtheorem{example}[theorem]{Example}

\newtheorem{proposition}[theorem]{Proposition}
\newtheorem{remark}[theorem]{Remark}

\newenvironment{proof}[1][Proof]{\noindent\textbf{#1.} }{\ \rule{0.5em}{0.5em}}
\input{tcilatex}

\begin{document}

\title{Interpolation properties for some scales of approximation spaces}
\author{Cristina Antonescu}
\date{}
\maketitle

\begin{abstract}
We obtain a result concerning the stability under the interpolation with
functional parameter method for the approximation spaces of
Lorentz-Marcinkiewicz type and also for the approximation spaces generated
by symmetric norming functions of a certain type.
\end{abstract}

\section{Introduction and general notions}

The subject of this paper is the real interpolation in approximation spaces.

Many of the approximation scales used in approximation theory are modeled by
some sequence ideals such as $l_{p},$ Lorentz, Lorentz-Zygmund or
Lorentz-Marcinkiewicz.

A natural framework which includes most of these scales and gives a unified
model for their study is provided by symmetric norming functions.
Unfortunately the abstract shape of a symmetric norming function is too
general to allow us to extend the basic interpolation results which hold for 
$l_{p},$ Lorentz, Lorentz-Zygmund or Lorentz-Marcinkiewicz scales to this
generality.

More precisely, while in the particular cases of $l_{p},$ Lorentz etc.
scales the stability under interpolation by the real method holds (the
classical version for $l_{p}$ and Lorentz cases, respectively the extended
one with functional parameter for Lorentz-Zygmund and Lorentz-Marcinkiewicz
cases), nothing is known about the behavior under the real interpolation of
the scale generated by symmetric norming functions .

Results given below give a partial extension of the known classical
stability properties to the general framework given by a scale of
approximation spaces generated by a certain class of symmetric norming
functions. In fact we prove the following result (for the notation see below)

The main observation which leads to the results below is that the Boyd
functions, which are the essential ingredient of the interpolation method
with functional parameter, have a shape close enough to the symmetric
norming functions of $\Phi ^{\varepsilon }$ type to allow us to obtain
informations about the interpolation process on the approximation spaces of $%
\Phi ^{\varepsilon }$ type (for the notation see below).

Before continuing, let us fix the notation used throughout the paper. By $%
E,F $ we denote Banach spaces over $\Gamma ,$ where $\Gamma $ is the real or
the complex field. We let $L(E,F):=\left\{ T:E\rightarrow F\mid T\text{ is
linear and bounded}\right\} .$ By $\widehat{k}$ we denote the set of all
decreasing positive sequences $x:=\left( x_{n}\right) _{n}$ such that $%
x_{n}=0$ eventually$.$We denote by $l_{\infty }$ the set of all scalar
sequences, $\left( x_{n}\right) _{n\in \mathbb{N}},$ with the property $%
\left\Vert x\right\Vert _{\infty }:=\underset{n\in \mathbb{N}}{\sup }%
\left\vert x_{n}\right\vert <\infty ,$ and by $c_{0}$ the set of all scalar
sequences, $\left( x_{n}\right) _{n\in \mathbb{N}},$ with the property $%
\underset{n\rightarrow \infty }{\lim }\left\vert x_{n}\right\vert =0.$ For $%
0<p<\infty $ we let $l_{p}$ be the set of all scalar sequences $\left(
x_{n}\right) _{n\in \mathbb{N}}$ such that $\left\Vert x\right\Vert
_{p}:=\left( \overset{\infty }{\underset{n=1}{\sum }}\left\vert
x_{n}\right\vert ^{p}\right) ^{\frac{1}{p}}<\infty .$

\section{Approximation schemes}

\subsection{Generalities}

We start with the general scheme described by a \textbf{quasi-normed abelian
group}, i.e. an abelian group $G$ endowed with a non-negative valued
function $\left\Vert \cdot \right\Vert :G\rightarrow \mathbb{R}_{+}$
satisfying:

\begin{enumerate}
\item $\left\Vert g\right\Vert >0,$ for every $g\in G,$ $g\neq 0;$

\item $\left\Vert -g\right\Vert =\left\Vert g\right\Vert ,$ for every $g\in
G;$

\item there is a constant, depending on $G,$ $c\geq 1$ such that 
\begin{equation*}
\left\Vert g+f\right\Vert \leq c\left( \left\Vert g\right\Vert +\left\Vert
f\right\Vert \right) ,
\end{equation*}
for every $g,f\in G.$
\end{enumerate}

\begin{definition}
([3], [15])) An \textbf{approximation scheme} on a quasi-normed abelian
group $\left( G,\left\Vert \cdot \right\Vert \right) $ consists of a pair $%
\left( G,\left( G_{n}\right) _{n\in \mathbb{N}}\right) $, where $\left(
G_{n}\right) _{n\in \mathbb{N}}$ is a sequence of subsets of $X$ satisfying
the following conditions:
\end{definition}

\begin{enumerate}
\item $G_{0}=\left\{ 0\right\} ;$

\item $G_{n}\subset G_{n+1},$ for every $n\in \mathbb{N};$

\item $G_{n}+G_{m}\subset G_{n+m},$ for every $n,m\in \mathbb{N};$

\item $G_{n}-G_{m}\subset G_{n+m},$ for every $n,m\in \mathbb{N}.$
\end{enumerate}

Associated with an approximation scheme $\left( G,\left( G_{n}\right) _{n\in 
\mathbb{N}}\right) $ one introduces the notion of an \textbf{approximation
space (}([3], [15])). We will use a function \newline
$E:G\rightarrow l_{\infty },$ defined by 
\begin{equation*}
E\left( g\right) :=\left( E_{n}\left( g\right) \right) _{n},\text{ for every 
}g\in G,
\end{equation*}%
where, $E_{n}\left( g\right) $ is the best approximation\textbf{\ }of $g$ by
elements of $G_{n-1},$ i.e. 
\begin{equation*}
E_{n}\left( g\right) :=\inf \left\{ \left\Vert g-h\right\Vert :h\in
G_{n-1}\right\} .
\end{equation*}%
The's basic properties of $E$ are as follows:

\begin{enumerate}
\item $\left\Vert g\right\Vert =E_{1}\left( g\right) \geq E_{2}\left(
g\right) \geq ...\geq 0,$ for every $g\in G$ ([15]);

\item $E_{n}\left( -g\right) =E_{n}\left( g\right) ,$ for every $n\in 
\mathbb{N}$ and $g\in G$ ([15]);

\item $E_{n+m}\left( f+g\right) \leq c\left( E_{n}\left( f\right)
+E_{m}\left( g\right) \right) ,$ for every $n,m\in \mathbb{N}$ and all $%
f,g\in G$ ([15])$.$
\end{enumerate}

\subsection{Examples}

We shall describe some important examples.

\begin{enumerate}
\item (\textbf{Operator ideals) }Let us take $G:=L(E,F)$ and $%
G_{n}:=F_{n}(E,F),$ where 
\begin{equation*}
F_{n}(E,F):=\left\{ T_{n}\in L(E,F):\dim T_{n}\leq n\right\} .
\end{equation*}%
Then for any $T\in L(E,F),$ the number $E_{n}\left( T\right) $ is the $n$-th 
\textbf{approximation number} of $T,$ denoted by $a_{n}(T)$ ([11], [13],
[14]).

\item \textbf{(Sequence ideals) }Let us now take $G:=l_{p}(0<p\leq \infty )$
and 
\begin{equation*}
G_{n}:=f_{p}^{\left( n\right) }:=\left\{ x:=\left( x_{m}\right) _{m}\in
l_{p}:\mathit{card\ }\left\{ x_{m}:x_{m}\neq 0\right\} \leq n\right\} .
\end{equation*}%
Then for any $x\in l_{p},$ the number $E_{n}\left( x\right) $ is the $n$-th 
\textbf{approximation number} of $x,$ denoted by $a_{n}(x).$ Let us remark
that, if the sequence $x:=\left( x_{n}\right) _{n}\in l_{p}$ is ordered such
that $\left\vert x_{n}\right\vert \geq \left\vert x_{n+1}\right\vert ,$ for
every $n,$ then $a_{n}\left( x\right) =\left\vert x_{n}\right\vert $ ([9],
[13], [14]).
\end{enumerate}

\subsection{Boyd functions}

One way to obtain approximation spaces is to use Boyd functions. We start
with a brief review of the notion of Boyd functions and the related
Lorentz-Marcienkiewicz scale.

\begin{definition}
$\left( [2],\text{ }\left[ 10\right] ,\text{ }[14]\right) $ We denote by 
\textbf{B} the class of all functions \newline
$\varphi :\left( 0,\infty \right) \rightarrow \left( 0,\infty \right) $
which have the following properties:
\end{definition}

\begin{enumerate}
\item $\varphi $ is continuous$;$

\item $\varphi \left( 1\right) =1;$

\item $\overline{\varphi }\left( t\right) :=\underset{s>0}{\sup }\frac{%
\varphi \left( ts\right) }{\varphi \left( s\right) }<\infty $ for any $t>0.$
\end{enumerate}

A straightforwards consequence of this definiton is that%
\begin{equation*}
\varphi \left( st\right) \leq \varphi \left( s\right) \overline{\varphi }%
\left( t\right) ,\text{ for any }s,t\in \left( 0,\infty \right) .
\end{equation*}

\begin{definition}
$\left( [2],\text{ }\left[ 10\right] ,\text{ }[14]\right) $ Given $\varphi $
in \textbf{B, the }$\overline{\varphi }$ \textbf{function's Boyd indices }$%
\alpha _{\overline{\varphi }},\beta _{\overline{\varphi }}$ are defined by: 
\begin{equation*}
\alpha _{\overline{\varphi }}:=\underset{1<t<\infty }{\inf }\frac{\log 
\overline{\varphi }\left( t\right) }{\log t}=\underset{t\rightarrow \infty }{%
\lim }\frac{\log \overline{\varphi }\left( t\right) }{\log t},
\end{equation*}%
and 
\begin{equation*}
\beta _{\overline{\varphi }}:=\underset{0<t<1}{\sup }\frac{\log \overline{%
\varphi }\left( t\right) }{\log t}=\underset{t\rightarrow 0}{\lim }\frac{%
\log \overline{\varphi }\left( t\right) }{\log t}.
\end{equation*}
\end{definition}

The Boyd indices satisfy the following relation

\begin{equation*}
-\infty <\beta _{\overline{\varphi }}\leq \alpha _{\overline{\varphi }%
}<\infty .
\end{equation*}

For future reference we collect some of the basic properties of the Boyd
functions in the following.

\begin{proposition}
If $\varphi ,\varphi _{1},\varphi _{2}$ are in \textbf{B}\textit{\ }and $a$
is a real number, then $\varphi _{1}\varphi _{2},\frac{\varphi _{1}}{\varphi
_{2}}$ are also in $\mathbf{B}$ and
\end{proposition}

\begin{enumerate}
\item $\overline{\varphi _{1}\varphi _{2}}\leq \overline{\varphi _{1}}\cdot 
\overline{\varphi _{2}}$

\item $\overline{\frac{\varphi _{1}}{\varphi _{2}}}\left( t\right) \leq 
\frac{\overline{\varphi _{1}}\left( t\right) }{\overline{\varphi _{2}}\left(
t\right) };$

\item $\beta _{\overline{\varphi }}>0$ if and only if $\underset{%
t\rightarrow 0}{\lim }\,\overline{\varphi }\left( t\right) =0$
\end{enumerate}

\begin{definition}
$\left( \left[ 14\right] ,\text{ }\left[ 15\right] \right) $ Let $\left(
G,\left( G_{n}\right) _{n}\right) $ be an approximation scheme, $\varphi \in 
$\textbf{B} and $0<q<\infty .$ An \textbf{approximation space }\ \textbf{of
Lorentz-Marcinkiewicz type} is defined as follows:%
\begin{equation*}
G_{\varphi ,q}:=\left\{ g\in G:\overset{\infty }{\underset{n=1}{\sum }}\left[
\varphi \left( n\right) E_{n}(g)\right] ^{q}n^{-1}<\infty \right\} .
\end{equation*}%
We also define an functional $\left\Vert \cdot \right\Vert _{\varphi
,q}:G_{\varphi ,q}\rightarrow \mathbb{R}$ by%
\begin{equation*}
\left\Vert g\right\Vert _{\varphi ,q}:=\left( \overset{\infty }{\underset{n=1%
}{\sum }}\left[ \varphi \left( n\right) E_{n}\left( g\right) \right]
^{q}n^{-1}\right) ^{\frac{1}{q}}\text{ for every }g\in G_{\varphi ,q}.
\end{equation*}
\end{definition}

We remark that $\left( G_{\varphi ,q},\left\Vert \cdot \right\Vert _{\varphi
,q}\right) $ is a quasi-normed abelian group, and for the particular case $%
\varphi \left( t\right) :=t^{p},$ $0<p<\infty ,$ we obtain the definition of
the classical approximation spaces $G_{p,q},$ where%
\begin{equation*}
\left( G_{p,q}:=\left\{ g\in G:\left\Vert g\right\Vert _{p,q}:=\left( 
\overset{\infty }{\underset{n=1}{\sum }}\left[ n^{p}E_{n}(g)\right]
^{q}n^{-1}\right) ^{\frac{1}{q}}<\infty \right\} \right) .
\end{equation*}

Now if we take for $\left( G,\left( G_{n}\right) _{n}\right) $ the
particular case $\left( l_{\infty },\left( f_{\infty }^{\left( n\right)
}\right) _{n}\right) $ we obtain the well known quasi-normed sequence ideal

\begin{equation*}
l_{\varphi ,q}:=\left\{ x\in l_{\infty }:\left\Vert x\right\Vert _{\varphi
,q}:=\left( \overset{\infty }{\underset{n=1}{\sum }}\left[ \varphi \left(
n\right) a_{n}(x)\right] ^{q}n^{-1}\right) ^{\frac{1}{q}}<\infty \right\} .
\end{equation*}%
Likewise, if we take for $\left( G,\left( G_{n}\right) _{n}\right) $ the
particular case $\left( L(E,F),\left( F_{n}(E,F)\right) _{n}\right) $ we
obtain the well known quasi-normed operator ideal 
\begin{equation*}
L_{\varphi ,q}:=\underset{E,F\text{ Banach spaces}}{\bigcup }L_{\varphi
,q}(E,F),
\end{equation*}%
where

\begin{equation*}
L_{\varphi ,q}(E,F):=\left\{ T\in L(E,F):\left\Vert T\right\Vert _{\varphi
,q}:=\left( \overset{\infty }{\underset{n=1}{\sum }}\left[ \varphi \left(
n\right) a_{n}(T)\right] ^{q}n^{-1}\right) ^{\frac{1}{q}}<\infty \right\} .
\end{equation*}

We note that equivalently we can define $G_{\varphi ,q}$ as the set of those 
$g$ in $G$ for which the sequence $\left( E_{n}\left( g\right) \right) _{n}$
belongs to $l_{\varphi ,q}.$

\subsection{Symmetric norming functions}

Another way for constructing approximation spaces uses a symmetric norming
function. We recall the definitions involved here.

\begin{definition}
([14], [15]) A function $\Phi :\widehat{k}\rightarrow \mathbb{R}$ is called
a \textbf{symmetric norming function} if there the following conditions are
fulfilled:

\begin{enumerate}
\item $\Phi \left( x\right) >0$ whenever $x\neq 0;$

\item $\Phi \left( \alpha x\right) =\alpha \Phi \left( x\right) $ for every $%
\alpha >0$ and $x$ in $\widehat{k};$

\item $\Phi \left( x+y\right) \leq \Phi \left( x\right) +\Phi \left(
y\right) $ for every $x$ and $y$ in $\widehat{k};$

\item $\Phi \left( \left\{ 1,0,...\right\} \right) =1;$

\item If $\overset{m}{\underset{n=1}{\sum }}x_{n}\leq \overset{m}{\underset{%
n=1}{\sum }}y_{n}$ for some $x=\left( x_{n}\right) _{n}$ and $y=\left(
x_{n}\right) _{n}$ in $\widehat{k}$ and for every $m$ in $\mathbb{N}^{\ast
}, $ then $\Phi \left( x\right) \leq \Phi \left( y\right) .$
\end{enumerate}
\end{definition}

\begin{remark}
$\left[ 14\right] $ Notice that the above definition can be extended to the
space $l_{\infty }$ of all bounded sequences in the following way. If $\Phi :%
\widehat{k}\rightarrow R$ is a symmetric norming function and $x:=\left(
x_{n}\right) _{n}\in l_{\infty }$ we define 
\begin{equation*}
\Phi \left( x\right) :=\underset{n}{\lim }\Phi \left( \left\{ a_{1}\left(
x\right) ,...,a_{n}\left( x\right) ,0,0,...\right\} \right) ,
\end{equation*}%
when $\underset{n}{\sup }$ $\Phi \left( \left\{ a_{1}\left( x\right)
,...,a_{n}\left( x\right) ,0,0,...\right\} \right) <\infty .$
\end{remark}

The most important examples of symmetric norming functions are the \textit{%
extremal }symmetric norming functions $\Phi _{1}$ and $\Phi _{\infty }$.
These are defined as follows:

\begin{equation*}
\Phi _{1}\left( x\right) :=\overset{n_{x}}{\underset{n=1}{\sum }}x_{n}\text{
and }\Phi _{\infty }\left( x\right) :=\underset{n}{\max }\text{ }x_{n},\text{
}x\in \widehat{k}.
\end{equation*}

It is easily seen that%
\begin{equation*}
\Phi _{\infty }\left( x\right) \leq \Phi \left( x\right) \leq \Phi
_{1}\left( x\right)
\end{equation*}%
for any symmetric norming function $\Phi $ and any $x\in \widehat{k}.$

In our future considerations a significant place is taken by the symmetric
norming functions of a certain type. These are the so-called $\Phi
^{\varepsilon }$ function. Their definition is presented in the next
proposition.

\begin{proposition}
$\left( \left[ 13\right] \right) $ Let $\varepsilon :=\left( \varepsilon
_{n}\right) _{n}$ be a decreasing sequence of positive real numbers with $%
\varepsilon _{1}=1.$ The function $\Phi ^{^{\varepsilon }}:\widehat{k}%
\rightarrow \mathbb{R}$ defined by%
\begin{equation*}
\Phi ^{^{\varepsilon }}\left( x\right) =\overset{n_{x}}{\underset{m=1}{\sum }%
}\varepsilon _{m}x_{m},\text{ for every }x:=\left( x_{m}\right) _{m\in 
\mathbb{N}}\text{ in }\widehat{k},
\end{equation*}%
is a symmetric norming function. If the sequence $\varepsilon $ has in
addition the properties $\underset{m\rightarrow \infty }{\lim }\varepsilon
_{m}=0$ and $\underset{m=1}{\overset{\infty }{\sum }}\varepsilon _{m}=\infty 
$ then $\Phi ^{^{\varepsilon }}\nsim \Phi _{\infty },$ and $\Phi
^{^{\varepsilon }}\nsim \Phi _{1}$ ( $\Phi \sim \Psi $ means $\underset{x\in 
\widehat{k}}{\sup }\frac{\Phi \left( x\right) }{\Psi \left( x\right) }%
<\infty ,$ and $\underset{x\in \widehat{k}}{\sup }\frac{\Psi \left( x\right) 
}{\Phi \left( x\right) }<\infty ).$ Let $\Phi ^{^{\varepsilon }}$ be a
function like above and $1\leq p<\infty .$ The function $\Phi _{\left(
p\right) }^{^{\varepsilon }}:\widehat{k}\rightarrow \mathbb{R}$ defined by 
\begin{equation*}
\Phi _{\left( p\right) }^{^{\varepsilon }}\left( x\right) =\left( \Phi
^{^{\varepsilon }}\left( \left( x_{m}^{p}\right) _{m}\right) \right) ^{\frac{%
1}{p}},\text{ for every }x:=\left( x_{m}\right) _{m}\text{ in }\widehat{k},
\end{equation*}%
is a symmetric function.
\end{proposition}

\begin{definition}
$\left( \left[ 5\right] \right) $ Let $\left( G,\left( G_{n}\right)
_{n}\right) $ be an approximation scheme and $\Phi $ a symmetric norming
function.\ An \textbf{approximation space of }$\Phi $\textbf{\ type} is
defined as follows:%
\begin{equation*}
G_{\Phi }:=\left\{ g\in G:\Phi \left( \left( E_{n}(g)\right) _{n}\right)
<\infty \right\} .
\end{equation*}%
We also define an functional $\left\Vert \cdot \right\Vert _{\Phi }:G_{\Phi
}\rightarrow \mathbb{R}$ by%
\begin{equation*}
\left\Vert g\right\Vert _{\Phi }:=\Phi \left( \left( E_{n}(g)\right)
_{n}\right) \text{ for every }g\in G_{\Phi }.
\end{equation*}
\end{definition}

We remark that $\left( G_{\Phi },\left\Vert \cdot \right\Vert _{\Phi
}\right) $ is a quasi-normed abelian group and for the particular case $\Phi
:=\Phi _{\left( q\right) }^{\alpha }$, with $\alpha :=\left( n^{pq-\frac{1}{q%
}}\right) ,$ we obtain, again, the definition of the classical approximation
spaces $G_{p,q}$.

Now if we take for $\left( G,\left( G_{n}\right) _{n}\right) $ the
particular case $\left( l_{\infty },\left( f_{\infty }^{\left( n\right)
}\right) _{n}\right) $ we obtain the well known quasi-normed sequence ideal

\begin{equation*}
l_{\Phi }:=\left\{ x\in l_{\infty }:\left\Vert x\right\Vert _{\Phi }:=\Phi
\left( \left( a_{n}(x)\right) _{n}\right) <\infty \right\} .
\end{equation*}%
Likewise, if we take for $\left( G,\left( G_{n}\right) _{n}\right) $ the
particular case $\left( L(E,F),\left( F_{n}(E,F)\right) _{n}\right) $ we
obtain the well known quasi-normed operator ideal 
\begin{equation*}
L_{\Phi }:=\underset{E,F\text{ Banach spaces}}{\bigcup }L_{\Phi }(E,F),
\end{equation*}%
where

\begin{equation*}
L_{\Phi }(E,F):=\left\{ T\in L(E,F):\left\Vert T\right\Vert _{\Phi }:=\Phi
\left( \left( a_{n}(T)\right) _{n}\right) <\infty \right\} ,
\end{equation*}

Equivalently we can define $G_{\Phi }$ as the set of those $g\in G$ for
which the sequence $\left( E_{n}\left( g\right) \right) _{n}$ belongs to $%
l_{\Phi }.$

\section{Stability results}

\subsection{Real interpolation}

Before continuing let us recall some results on real interpolation.

We consider couples $\left( E_{0},E_{1}\right) $ of quasi-normed spaces $%
E_{0}$ and $E_{1}$, which are both continuously embedded in a quasi-normed
space $E.$ This means that $E_{i}\subset E$ and there is a constant $c_{i}$
such that $\left\Vert x\right\Vert _{E}\leq c_{i}\left\Vert x\right\Vert
_{E_{i}}$ for every $x$ in $E_{i},$ $i\in \left\{ 0,1\right\} .$ In the
sequel we let $\hookrightarrow $ denote a continuous embedding. We say that
such a couple $\left( E_{0},E_{1}\right) $ is a \textbf{quasi-normed
interpolation couple}.

Let $\left( E_{0},E_{1}\right) $ be a quasi-normed interpolating couple. For
every $x$ in $E_{0}+E_{1}$ J. Peetre defined the functional 
\begin{equation*}
K\left( t,x,E_{0},E_{1}\right) =K(t,x):=\underset{x=x_{0}+x_{1}}{\inf }%
\left( \left\Vert x_{0}\right\Vert _{E_{0}}+t\left\Vert x_{1}\right\Vert
_{E_{1}}\right) ,
\end{equation*}%
where $x_{i}\in E_{i},$ $i\in \left\{ 0,1\right\} $ and $0<t<\infty .$ Let
now $\left( E_{0},E_{1}\right) $ be a quasi-normed interpolation couple, $%
0<q<\infty $ and $\varphi \in \mathbf{B}.$ We shall consider the set

\begin{equation*}
\left( E_{0},E_{1}\right) _{\varphi ,q}:=\left\{ x\in E_{0}+E_{1}:\overset{%
\infty }{\underset{0}{\int }}\left[ \varphi \left( t\right) ^{-1}K\left(
t,x\right) \right] ^{q}\frac{dt}{t}<\infty \right\} .
\end{equation*}%
It is important to notice that for any $\varphi \in $\textbf{B} and $%
0<q<\infty ,$ the functional $\left\Vert \cdot \right\Vert _{\varphi
,q}:\left( E_{0},E_{1}\right) _{\varphi ,q}\rightarrow \mathbb{R}_{+},$
defined by

\begin{equation*}
\left\Vert x\right\Vert _{\left( E_{0},E_{1}\right) _{\varphi ,q}}:=\left( 
\overset{\infty }{\underset{0}{\int }}\left[ \varphi \left( t\right)
^{-1}K\left( t,x\right) \right] ^{q}\frac{dt}{t}\right) ^{\frac{1}{q}},
\end{equation*}%
for every $x\in \left( E_{0},E_{1}\right) _{\varphi ,q},$ is a quasi-norm$%
\mathbf{.}$

Consider now the space $E_{\Sigma }:=E_{0}+E_{1}$ equipped with the
quasi-norm 
\begin{equation*}
\left\Vert x\right\Vert _{\Sigma }:=\underset{%
\begin{array}{c}
x=x_{0}+x_{1}, \\ 
x_{i}\in E_{i}%
\end{array}%
}{\inf }\left( \left\Vert x_{0}\right\Vert _{0}+\left\Vert x_{1}\right\Vert
_{1}\right) =K(1,x)
\end{equation*}%
and also the space $E_{\Delta }:=E_{0}\dbigcap E_{1}$ equipped with the
quasi-norm 
\begin{equation*}
\left\Vert x\right\Vert _{\Delta }:=\max \left( \left\Vert x\right\Vert
_{0},\left\Vert x\right\Vert _{1}\right) .
\end{equation*}%
We remark that 
\begin{equation}
\left( E_{0},E_{1}\right) _{\varphi ,q}\hookrightarrow E_{\Sigma }
\end{equation}%
and, 
\begin{equation}
\text{if }0<\beta _{\overline{\varphi }}\leq \alpha _{\overline{\varphi }}<1%
\text{ then }E_{\Delta }\hookrightarrow \left( E_{0},E_{1}\right) _{\varphi
,q}.
\end{equation}

The above construction which starts with an interpolation couple $\left(
E_{0},E_{1}\right) $ and give us the space $\left( E_{0},E_{1}\right)
_{\varphi ,q}$ (called \textbf{real interpolation method with functional
parameter) }was introduced by T.F. Kalugina ([2], [4], [6], [10]) as an
extension of the classical real method due to J. Peetre ([1], [11], [13],
[14]). More precisely, if we take $\varphi (t)=t^{-\theta },$ where $%
0<\theta <1$, we get the classical real interpolation space $\left(
E_{0},E_{1}\right) _{\theta ,q}.$

We recall now the \textbf{reiteration theorem} of the real interpolation
method with functional parameter which is the main ingredient of our proofs.

\begin{theorem}
$\left( \left[ 10\right] \right) $ Let $\overline{A}=\left\{
A_{0},A_{1}\right\} $ be an interpolation couple of quasi-normed spaces.
Take $f,f_{0},f_{1}$ in \textbf{B}\textit{, \ where }$f$ in addition$\ $%
satisfies $0<\beta _{\overline{f}}\leq \alpha _{\overline{f}}<1$ and let $u=%
\frac{f_{1}}{f_{0}},$ $g=f_{0}\left( f\circ u\right) ,$ $0<p,q_{0},q_{1}\leq
\infty .$ By $E_{i}$ we denote the interpolation space $\left(
A_{0},A_{1}\right) _{f_{i},q_{i}},$ where $i\in \left\{ 0,1\right\} .$ If
one of the following two hypotheses is fulfilled:\newline
1. $\beta _{\overline{f_{0}}}>0$ in the case $q_{0}<\infty ,$ respectively $%
\underset{t\leq 1}{\sup }f_{0}<\infty $ in the case $q_{0}=\infty ,$ and $%
\alpha _{\overline{f_{1}}}<1$ in the case $q_{1}<\infty ,$ respectively $%
\underset{t\geq 1}{\sup }\frac{f_{1}\left( t\right) }{t}<\infty $ in the
case $q_{1}=\infty ,$ when $\beta _{\overline{u}}>0$ or \newline
2. $\beta _{\overline{f_{1}}}>0$ in the case $q_{1}<\infty ,$ respectively $%
\underset{t\leq 1}{\sup }f_{1}<\infty $ in the case $q_{1}=\infty ,$ and $%
\alpha _{\overline{f_{0}}}>0$ in the case $q_{0}<\infty ,$ respectively $%
\underset{t\geq 1}{\sup }\frac{f_{0}\left( t\right) }{t}<\infty $ in the
case $q_{0}=\infty ,$ when $\alpha _{\overline{u}}<0,$ \newline
then 
\begin{equation*}
g\in \mathbf{B\ }\text{and }\left( E_{0},E_{1}\right) _{f,p}=\left(
A_{0},A_{1}\right) _{g,p}.
\end{equation*}
\end{theorem}

Let us mention that, from this abstract reiteration theorem, F. Cobos has
obtained the stability under functional parameter interpolation's process
for the operator ideals introduced by him.

\begin{theorem}
$\left( \left[ 2\right] \right) $ Take $E,F$ Banach spaces, the numbers $%
q_{0},q_{1},q$ in $(0,\infty ]$ and the functions $\chi ,\varphi
_{0},\varphi _{1}$ in $\mathbf{B}$. Let now consider the functions $\varphi
:\left( 0,\infty \right) \rightarrow \left( 0,\infty \right) ,$ $\rho
:\left( 0,\infty \right) \rightarrow \left( 0,\infty \right) $ defined by 
\begin{equation*}
\varphi \left( t\right) :=\frac{\varphi _{0}\left( t\right) }{\varphi
_{1}\left( t\right) },
\end{equation*}%
respectively by 
\begin{equation*}
\rho \left( t\right) =\frac{\varphi _{0}\left( t\right) }{\chi \left(
\varphi \left( t\right) \right) }.
\end{equation*}%
If $0<\beta _{\overline{\chi }}\leq $ $\alpha _{\overline{\chi }}<1,\beta _{%
\overline{\varphi _{i}}}>0\left( i=0,1\right) $ and $\beta _{\overline{%
\varphi }}>0$ or $\alpha _{\overline{\varphi }}<0,$ then 
\begin{equation*}
\rho \in \mathbf{B}\text{ and }\left( L_{\varphi _{0},q_{0}}(E,F),L_{\varphi
_{1},q_{1}}(E,F)\right) _{\chi ,q}=L_{\rho ,q}(E,F),
\end{equation*}%
with equivalent quasi-norms.
\end{theorem}

\subsection{Stability results for approximation spaces: The
Lorentz-Marcinkiewicz case}

The main result of this section is as follows.

\begin{theorem}
Take the numbers $q_{0},q_{1},q$ in $(0,\infty ]$ and the functions $\chi
,\varphi _{0},\varphi _{1}$ in $\mathbf{B}$. Let now consider the functions $%
\varphi ,$ $\rho :\left( 0,\infty \right) \rightarrow \left( 0,\infty
\right) $ defined by 
\begin{equation*}
\varphi \left( t\right) :=\frac{\varphi _{0}\left( t\right) }{\varphi
_{1}\left( t\right) },\text{ respectively }\rho \left( t\right) :=\frac{%
\varphi _{0}\left( t\right) }{\chi \left( \varphi \left( t\right) \right) }
\end{equation*}
If the following conditions are fulfilled: 
\begin{equation}
0<\beta _{\overline{\chi }}\leq \alpha _{\overline{\chi }}<1,\beta _{%
\overline{\varphi _{i}}}>0\text{, where }i\in \left\{ 0,1\right\}
\end{equation}
and 
\begin{equation}
\beta _{\overline{\varphi }}>0\text{ or }\alpha _{\overline{\varphi }}<0
\end{equation}
then%
\begin{equation*}
\rho \in \mathbf{B}\text{ and }\left( G_{\varphi _{0},q_{0}},G_{\varphi
_{1},q_{1}}\right) _{\chi ,q}=G_{\rho ,q},
\end{equation*}%
with equivalent quasi-norms.
\end{theorem}

The proof is constructed in two steps, the first of which has independent
interest. We shall use the notation%
\begin{equation*}
G_{p}:=\left\{ g\in G:\overset{\infty }{\underset{n=1}{\sum }}\left[
E_{n}\left( g\right) \right] ^{p}<\infty \right\} .
\end{equation*}

\begin{theorem}
Take $\varphi $ in $\mathbf{B}$ which satisfies the condition $0<\beta _{%
\overline{\varphi }}$ and $q$ in $(0,\infty ].$ If $0<p_{0}<p_{1}\leq \infty 
$ are such that 
\begin{equation*}
\frac{1}{p_{1}}<\beta _{\overline{\varphi }}<\alpha _{\overline{\varphi }}<%
\frac{1}{p_{0}}
\end{equation*}
and if $\rho :\left( 0,\infty \right) \rightarrow \left( 0,\infty \right) $
is defined by 
\begin{equation*}
\rho \left( t\right) =t^{\frac{p_{1}}{p_{0}-p_{1}}}\left( \varphi \left( t^{%
\frac{p_{0}p_{1}}{p_{1}-p_{0}}}\right) \right) ^{-1}
\end{equation*}%
(in the case $p_{1}<\infty ),$ respectively 
\begin{equation*}
\rho \left( t\right) =t\left( \varphi \left( t^{p_{0}}\right) \right) ^{-1}
\end{equation*}%
(in the case $p_{1}=\infty ),$ then 
\begin{equation*}
\rho \in \mathbf{B}\text{ and }\left( G_{p_{0}},G_{p_{1}}\right) _{\rho
,q}=G_{\varphi ,q},
\end{equation*}%
with equivalent quasi-norms.
\end{theorem}

\begin{proof}
The fundamental observation is that the equivalence 
\begin{equation*}
K\left( t,T,L_{p_{0}}^{\left( a\right) }(E,F),L_{p_{1}}^{\left( a\right)
}(E,F)\right) \simeq K\left( t,\left( a_{n}\left( T\right) \right)
_{n},l_{p_{0}},l_{p_{1}}\right) ,
\end{equation*}%
proved by H. K\"{o}nig in [7], where \textquotedblright $\simeq "$ indicates
equivalence with constants that do not depend on $t$ or $T,$ remain valid
for the abstract case of the application $E:X\rightarrow l_{\infty },$ the
proof being the same. Hence 
\begin{equation*}
K\left( t,g,G_{p_{0}},G_{p_{1}}\right) \simeq K\left( t,\left( E_{n}\left(
g\right) \right) _{n},l_{p_{0}},l_{p_{1}}\right) .
\end{equation*}%
Now we obtain the following equivalences%
\begin{equation*}
g\in \left( G_{p_{0}},G_{p_{1}}\right) _{\rho ,q}\Leftrightarrow \overset{%
\infty }{\underset{0}{\int }}\left[ \rho \left( t\right) ^{-1}K\left(
t,f,G_{p_{0}},G_{p_{1}}\right) \right] ^{q}\frac{dt}{t}<\infty
\Leftrightarrow
\end{equation*}%
\begin{equation*}
\Leftrightarrow \overset{\infty }{\underset{0}{\int }}\left[ \rho \left(
t\right) ^{-1}K\left( t,\left( E_{n}\left( f\right) \right)
_{n},l_{p_{0}},l_{p_{1}}\right) \right] ^{q}\frac{dt}{t}<\infty
\Leftrightarrow \left( E_{n}\left( f\right) \right) _{n}\in \left(
l_{p_{0}},l_{p_{1}}\right) _{\rho ,q}\Leftrightarrow
\end{equation*}%
\begin{equation*}
\Leftrightarrow \left( E_{n}\left( g\right) \right) _{n}\in l_{\varphi
,q}\Leftrightarrow \underset{n}{\sum }\left[ \varphi \left( n\right)
E_{n}\left( g\right) \right] ^{q}n^{-1}<\infty \Leftrightarrow g\in
G_{\varphi ,q}.
\end{equation*}
\end{proof}

From the above result, applying Theorem 10, we can prove the main result of
this section

\begin{proof}
(of Theorem 12) From the previous theorem we conclude that 
\begin{equation*}
G_{\varphi _{i},q_{i}}=\left( G_{r},G\right) _{f_{i},q_{i}},\text{ for }i\in
\left\{ 0,1\right\} .
\end{equation*}%
On the other hand we see that the hypotheses of Theorem 10\textbf{\ }are
fulfilled. In our case%
\begin{equation*}
f=\chi ,f_{i}\left( t\right) =\frac{t}{\varphi _{i}\left( t^{r}\right) }%
(i\in \left\{ 0,1\right\} ),
\end{equation*}%
the number $r$ being chosen in such a way that $\frac{1}{r}>\max \left\{
\alpha _{\overline{\varphi _{0}}},\alpha _{\overline{\varphi _{1}}},\alpha _{%
\overline{\rho }}\right\} ,$ and 
\begin{equation*}
u\left( t\right) =\frac{f_{1}\left( t\right) }{f_{0}\left( t\right) }%
=\varphi \left( t^{r}\right) .
\end{equation*}%
Applying now Theorem 10 we obtain%
\begin{equation*}
\left( G_{\varphi _{0},q_{0}},G_{\varphi _{1},q_{1}}\right) _{\chi
,q}=\left( \left( G_{r},G\right) _{f_{0},q_{0}},\left( G_{r},G\right)
_{f_{1},q_{1}}\right) _{\chi ,q}=\left( G_{r},G\right) _{g,q},
\end{equation*}%
where 
\begin{equation*}
g\left( t\right) =f_{0}\left( t\right) \left( f\circ u\right) \left(
t\right) =\frac{t}{\varphi _{0}\left( t^{r}\right) }\chi \left( \varphi
\left( t^{r}\right) \right) .
\end{equation*}%
Applying again the previous theorem, this time with $\rho =g$, we obtain 
\begin{equation*}
\left( G_{r},G\right) _{g,q}=G_{\varpi ,q},
\end{equation*}%
where%
\begin{equation*}
\varpi \left( t^{r}\right) =\frac{t}{g\left( t\right) }=\frac{\varphi
_{0}\left( t^{r}\right) }{\chi \left( \varphi \left( t^{r}\right) \right) }.
\end{equation*}%
Hence we have obtained 
\begin{equation*}
\varpi \left( t\right) =\frac{\varphi _{0}\left( t\right) }{\chi \left(
\varphi \left( t\right) \right) }=\rho \left( t\right) .
\end{equation*}%
In conclusion,%
\begin{equation*}
\left( G_{r},G\right) _{g,q}=G_{\rho ,q},
\end{equation*}%
and moreover we can write the desired equality%
\begin{equation*}
\left( G_{\varphi _{0},q_{0}},G_{\varphi _{1},q_{1}}\right) _{\chi
,q}=G_{\rho ,q}.
\end{equation*}
\end{proof}

\subsubsection{Stability results for approximation schemes: The symmetric
norming functions case}

Speaking now about the frame constructed with symmetric norming functions
the basic idea for finding interpolation results for the approximation
spaces of of type $G_{\Phi }$ is that the symmetric norming functions of $%
\Phi ^{\varepsilon }$ can be arranged as some Boyd functions like we shall
present in the following construction.

\begin{definition}
Let $\alpha =\left( \alpha _{n}\right) _{n}$ be a sequence of real numbers
with the following properties:\newline
1. $1=\alpha _{1}\geq \alpha _{2}\geq ...\geq 0;\newline
2.$ $\underset{n\rightarrow \infty }{\lim }\alpha _{n}=0;\newline
3.$ $\overset{\infty }{\underset{n=1}{\sum }}\alpha _{n}=\infty ;\newline
4.$ $M\left( p\right) :=\underset{n}{\sup }\frac{\alpha _{\left[ \frac{n}{p}%
\right] }}{\alpha _{n}}<\infty $ for every fixed $p>1,$where by $\left[ t%
\right] $ we denote the greatest integer less or equal than $t.\newline
$For every sequence like above and every positive number $p$ we define the
function $\varphi _{\alpha ,p}:\left( 0,\infty \right) \rightarrow \left(
0,\infty \right) $ as follows 
\begin{equation*}
\varphi _{\alpha ,p}\left( t\right) :=\left\{ 
\begin{array}{l}
t^{\frac{1}{p}}\text{\qquad if }t\in \left( 0,1\right) \\ 
\left( \alpha _{t}t\right) ^{\frac{1}{p}}\text{ if }t\in \mathbb{N}^{\ast }
\\ 
\left( 1-\left\{ t\right\} \right) \varphi _{\alpha ,p}\left( \left[ t\right]
\right) +\left\{ t\right\} \varphi _{\alpha ,p}\left( \left[ t\right]
+1\right) \text{ if }t\in \left( 1,\infty \right) \backslash \mathbb{N}%
\end{array}%
\right.
\end{equation*}%
where $\left\{ t\right\} :=t-\left[ t\right] .$
\end{definition}

To prove that Definition 14 is not void we present the following.

\begin{example}
The sequence $\left( \frac{1}{n^{a}}\right) _{n},$ $a\leq 1$ has the
properties 1-4 from Definition 14.
\end{example}

\begin{proof}
The properties 1-3 are obviously fulfilled. To verify the fourth condition
fix $p$ in $\mathbb{N},$ $p>1.$ It is known that%
\begin{equation*}
\left[ \frac{s}{p}\right] =\left[ \frac{\left[ s\right] }{p}\right] =\left[ %
\left[ s\right] \frac{1}{p}\right] \geq \left[ s\right] \left[ \frac{1}{p}%
\right]
\end{equation*}%
for every positive number $s.$ Hence we shall obtain 
\begin{equation*}
\alpha _{\left[ \frac{s}{p}\right] }\leq \alpha _{\left[ s\right] \left[ 
\frac{1}{p}\right] }
\end{equation*}
and furthermore 
\begin{equation*}
\frac{\alpha _{\left[ \frac{n}{p}\right] }}{\alpha _{n}}\leq \frac{\alpha _{n%
\left[ \frac{1}{p}\right] }}{\alpha _{n}}=\frac{1}{\left[ \frac{1}{p}\right]
^{\alpha }}.
\end{equation*}%
In conclusion 
\begin{equation*}
\underset{n}{\sup }\,\frac{\alpha _{\left[ \frac{n}{p}\right] }}{\alpha _{n}}%
\leq \frac{1}{\left[ \frac{1}{p}\right] ^{\alpha }}<\infty .
\end{equation*}
\end{proof}

It is obvious that the function $\Phi ^{^{\alpha }}:l_{\infty }\rightarrow 
\mathbb{R}$ defined by

\begin{equation*}
\Phi ^{^{\alpha }}\left( x\right) :=\overset{\infty }{\underset{n=1}{\sum }}%
\alpha _{n}a_{n}\left( x\right) ,
\end{equation*}%
for all $x\in l_{\infty }$ is a symmetric norming function. More interesting
is the fact that the above definition leads to a Boyd function.

\begin{theorem}
Let $\alpha =\left( \alpha _{n}\right) _{n}$ be a sequence which has the
properties 1-4 from Definition 14 and $p$ $\in \left( 0,\infty \right) .$
Then $\varphi _{\alpha ,p}\in \mathbf{B}.$
\end{theorem}

\begin{proof}
We shall verify the axioms from the definition of a \textbf{B}-function. The
construction ensures the continuity of the function $\varphi _{\alpha ,p}$
and also the equality 
\begin{equation*}
\varphi _{\alpha ,p}\left( 1\right) =\left( \alpha _{1}\right) ^{\frac{1}{p}%
}=1.
\end{equation*}%
In order verify the third condition we consider an arbitrary number $t\in
\left( 0,\infty \right) $. We shall evaluate from above $\frac{\varphi
_{\alpha ,p}\left( st\right) }{\varphi _{\alpha ,p}\left( s\right) }$ for $%
s\in \left( 0,\infty \right) .$ This is way in the sequel shall take account
of the inequality $\left[ xy\right] \geq \left[ x\right] \left[ y\right] $
which is true for all positive $x,y.$ Without loss of generality we may
assume that $t\neq 1.$ There are two cases. The first one is to consider all 
$s\in \left( 0,\infty \right) $ such that $st<1.$ Then we'll have 
\begin{equation*}
\frac{\varphi _{\alpha ,p}\left( st\right) }{\varphi _{\alpha ,p}\left(
s\right) }=\frac{\left( st\right) ^{\frac{1}{p}}}{s^{\frac{1}{p}}}=t^{\frac{1%
}{p}}.
\end{equation*}%
The second one case is to consider all $s\in \left( 0,\infty \right) $ such
that $st>1.$ In that case we have to analyze two situations. For the moment
we fix $t<1$. Then we obtain the following relations:%
\begin{equation*}
\frac{\varphi _{\alpha ,p}\left( st\right) }{\varphi _{\alpha ,p}\left(
s\right) }=\frac{\left( 1-\left\{ st\right\} \right) \varphi _{\alpha
,p}\left( \left[ st\right] \right) +\left\{ st\right\} \varphi _{\alpha
,p}\left( \left[ st\right] +1\right) }{\left( 1-\left\{ s\right\} \right)
\varphi _{\alpha ,p}\left( \left[ s\right] \right) +\left\{ s\right\}
\varphi _{\alpha ,p}\left( \left[ s\right] +1\right) }=
\end{equation*}%
\begin{equation*}
=\frac{\left( 1-\left\{ st\right\} \right) \left( \alpha _{\left[ st\right]
}\left( \left[ st\right] \right) \right) ^{\frac{1}{p}}+\left\{ st\right\}
\left( \alpha _{\left[ st\right] +1}\left( \left[ st\right] +1\right)
\right) ^{\frac{1}{p}}}{\left( 1-\left\{ s\right\} \right) \left( \alpha _{%
\left[ s\right] }\left[ s\right] \right) ^{\frac{1}{p}}+\left\{ s\right\}
\left( \alpha _{\left[ s\right] +1}\left( \left[ s\right] +1\right) \right)
^{\frac{1}{p}}}=
\end{equation*}%
\begin{equation*}
=\frac{\left( 1-\left\{ st\right\} \right) \left( \alpha _{\left[ st\right]
}\left( \left[ st\right] \right) \right) ^{\frac{1}{p}}+\left\{ st\right\}
\left( \alpha _{\left[ st\right] +1}\left( \left[ st\right] +1\right)
\right) ^{\frac{1}{p}}}{\left( 1-\left\{ s\right\} \right) \left( \alpha _{%
\left[ s\right] }\left[ st\frac{1}{t}\right] \right) ^{\frac{1}{p}}+\left\{
s\right\} \left( \alpha _{\left[ s\right] +1}\left( \left[ st\frac{1}{t}%
\right] +1\right) \right) ^{\frac{1}{p}}}\leq
\end{equation*}%
\begin{equation*}
\leq \frac{\alpha _{\left[ st\right] }^{\frac{1}{p}}\left( \left( 1-\left\{
st\right\} \right) \left[ st\right] ^{\frac{1}{p}}+\left\{ st\right\} \left( %
\left[ st\right] +1\right) ^{\frac{1}{p}}\right) }{\alpha _{\left[ s\right]
+1}^{\frac{1}{p}}\left( \left( 1-\left\{ s\right\} \right) \left[ st\right]
^{\frac{1}{p}}\left[ \frac{1}{t}\right] ^{\frac{1}{p}}+\left\{ s\right\}
\left( \left[ st\right] \left[ \frac{1}{t}\right] +1\right) ^{\frac{1}{p}%
}\right) }=
\end{equation*}%
\begin{equation*}
=\left( \frac{\alpha _{\left[ s\cdot t\right] }}{\alpha _{\left[ s\right] +1}%
}\right) ^{\frac{1}{p}}\frac{\left( 1-\left\{ st\right\} \right) \left[ st%
\right] ^{\frac{1}{p}}+\left\{ st\right\} \left( \left[ st\right] +1\right)
^{\frac{1}{p}}}{\left( 1-\left\{ s\right\} \right) \left[ st\right] ^{\frac{1%
}{p}}\left[ \frac{1}{t}\right] ^{\frac{1}{p}}+\left\{ s\right\} \left( \left[
st\right] \left[ \frac{1}{t}\right] +1\right) ^{\frac{1}{p}}}=
\end{equation*}%
\begin{equation*}
=\left( \frac{\alpha _{\left[ st\right] }}{\alpha _{\left[ s\right] +1}}%
\right) ^{\frac{1}{p}}\frac{1}{\left[ \frac{1}{t}\right] ^{\frac{1}{p}}}%
\cdot \frac{\left( 1-\left\{ st\right\} \right) +\left\{ st\right\} \left( 1+%
\frac{1}{\left[ st\right] }\right) ^{\frac{1}{p}}}{\left( 1-\left\{
s\right\} \right) +\left\{ s\right\} \left( 1+\frac{1}{\left[ st\right] %
\left[ \frac{1}{t}\right] }\right) ^{\frac{1}{p}}}\leq
\end{equation*}%
\begin{equation*}
\leq \left[ M\left( \frac{1}{t}\right) \right] ^{\frac{1}{p}}\frac{1}{\left[ 
\frac{1}{t}\right] ^{\frac{1}{p}}}\frac{\left( 1-\left\{ st\right\} \right)
+\left\{ st\right\} \left( 1+\frac{1}{\left[ st\right] }\right) ^{\frac{1}{p}%
}}{\left( 1-\left\{ s\right\} \right) +\left\{ s\right\} \left( 1+\frac{1}{%
\left[ st\right] \left[ \frac{1}{t}\right] }\right) ^{\frac{1}{p}}}.
\end{equation*}%
Now we assume that $t>1.$ Then%
\begin{equation*}
\frac{\varphi _{\alpha ,p}\left( st\right) }{\varphi _{\alpha ,p}\left(
s\right) }=\frac{\left( 1-\left\{ st\right\} \right) \varphi _{\alpha
,p}\left( \left[ st\right] \right) +\left\{ st\right\} \varphi _{\alpha
,p}\left( \left[ st\right] +1\right) }{\left( 1-\left\{ s\right\} \right)
\varphi _{\alpha ,p}\left( \left[ s\right] \right) +\left\{ s\right\}
\varphi _{\alpha ,p}\left( \left[ s\right] +1\right) }=
\end{equation*}%
\begin{equation*}
=\frac{\left( 1-\left\{ st\right\} \right) \left( \alpha _{\left[ st\right]
}\left( \left[ st\right] \right) \right) ^{\frac{1}{p}}+\left\{ st\right\}
\left( \alpha _{\left[ st\right] +1}\left( \left[ st\right] +1\right)
\right) ^{\frac{1}{p}}}{\left( 1-\left\{ s\right\} \right) \left( \alpha _{%
\left[ s\right] }\left[ s\right] \right) ^{\frac{1}{p}}+\left\{ s\right\}
\left( \alpha _{\left[ s\right] +1}\left( \left[ s\right] +1\right) \right)
^{\frac{1}{p}}}\leq
\end{equation*}%
\begin{equation*}
\leq \frac{\alpha _{\left[ s\right] }^{\frac{1}{p}}\left( \left[ st\right]
+1\right) ^{\frac{1}{p}}}{\alpha _{\left[ s\right] +1}^{\frac{1}{p}}\left[ s%
\right] ^{\frac{1}{p}}}=\left( \frac{\alpha _{\left[ s\right] }}{\alpha _{%
\left[ s\right] +1}}\right) ^{\frac{1}{p}}\left( \frac{\left( \left[ st%
\right] +1\right) }{\left[ s\right] }\right) ^{\frac{1}{p}}\leq
\end{equation*}%
\begin{equation*}
\leq \left( \frac{\left( \left[ st\right] +1\right) }{\left[ s\right] }%
\right) ^{\frac{1}{p}}\leq \left( \frac{st+1}{s-1}\right) ^{\frac{1}{p}}.
\end{equation*}%
In conclusion 
\begin{equation*}
\underset{s>0}{\sup }\,\frac{\varphi _{\alpha ,p}\left( st\right) }{\varphi
_{\alpha ,p}\left( s\right) }<\infty \text{ for every }t\in \left( 0,\infty
\right) .
\end{equation*}
\end{proof}

We are prepared now for the main result of this section.

\begin{theorem}
Consider an approximation scheme $\left( G,\left( G_{n}\right) _{n\in 
\mathbb{N}}\right) $. Let $\alpha :=\left( \alpha _{n}\right) _{n}$ and $%
\beta :=\left( \beta _{n}\right) _{n}$ be sequences having the properties
1-4 from Definition 14. If $\alpha ,\beta $ in addition satisfy the
conditions\newline
1. $\beta _{n}\leq \alpha _{n},$ for every $n\in \mathbb{N}^{\ast }$ and $%
\newline
2.$ $\underset{t\rightarrow 0}{\lim }M_{\alpha }\left( \frac{1}{t}\right)
t=0,$ respectively $\underset{t\rightarrow 0}{\lim }M_{\beta }\left( \frac{1%
}{t}\right) t=0$ $\newline
$and we choose the positive numbers $p,q,l$ satisfying the following four
relations 
\begin{equation*}
1\leq p\leq q<\infty ,\text{ }l>1,\text{ }p+ql>q\text{ and }\frac{pql}{p+ql-q%
}>1,
\end{equation*}%
then 
\begin{equation*}
\left( G_{\Phi _{\left( p\right) }^{^{^{\alpha }}}},G_{\Phi _{\left(
q\right) }^{^{^{\beta }}}}\right) _{f,r}=G_{\Phi _{\left( r\right)
}^{^{^{\gamma }}}},\text{ with equivalent quasi-norms,}
\end{equation*}%
where $f:\left( 0,\infty \right) \rightarrow \left( 0,\infty \right) $ is
given by 
\begin{equation*}
f\left( t\right) :=t^{\frac{1}{l}}
\end{equation*}%
and $\gamma :=\left( \gamma _{n}\right) _{n}$ is given by 
\begin{equation*}
\gamma _{n}:=\alpha _{n}^{r\left( \frac{1}{p}-\frac{1}{pl}\right) }\beta
_{n}^{\frac{r}{ql}}.
\end{equation*}
\end{theorem}

\begin{proof}
We start by proving that $l_{\varphi _{\alpha ,p},p}=l_{\Phi _{\left(
p\right) }^{^{\alpha }}}.$ From the definitions 
\begin{equation*}
l_{\Phi _{\left( p\right) }^{^{\alpha }}}=\left\{ x\in l_{\infty }:\left( 
\overset{\infty }{\underset{n=1}{\sum }}\alpha _{n}\left[ a_{n}\left(
x\right) \right] ^{p}\right) ^{\frac{1}{p}}<\infty \right\}
\end{equation*}%
and 
\begin{equation*}
l_{\varphi _{\alpha ,p},p}=\left\{ x\in l_{\infty }:\left( \overset{\infty }{%
\underset{n=1}{\sum }}\left[ \varphi _{\alpha ,p}\left( n\right) a_{n}\left(
x\right) \right] ^{p}n^{-1}\right) ^{\frac{1}{p}}<\infty \right\} .
\end{equation*}%
Consequently we obtain%
\begin{equation*}
l_{\varphi _{\alpha ,p},p}=\left\{ x\in l_{\infty }:\left( \overset{\infty }{%
\underset{n=1}{\sum }}\left[ \left( \alpha _{n}n\right) ^{\frac{1}{p}%
}a_{n}\left( x\right) \right] ^{p}n^{-1}\right) ^{\frac{1}{p}}<\infty
\right\} =
\end{equation*}%
\begin{equation*}
=\left\{ x\in l_{\infty }:\left( \overset{\infty }{\underset{n=1}{\sum }}%
\alpha _{n}\left[ a_{n}\left( x\right) \right] ^{p}\right) ^{\frac{1}{p}%
}<\infty \right\} =l_{\Phi _{\left( p\right) }^{^{\alpha }}}.
\end{equation*}%
It is easy to check $\left\Vert \cdot \right\Vert _{\Phi _{\left( p\right)
}^{^{\alpha }}}=\left\Vert \cdot \right\Vert _{\varphi _{\alpha ,p},p}.$
Similarly $l_{\Phi _{\left( q\right) }^{^{\beta }}}=l_{\varphi _{\beta
,q},q}.$\newline
Let $g\in G.$ From the definition of the approximation spaces $G_{\Phi
_{\left( p\right) }^{^{^{\alpha }}}}$ we obtain the following equivalences%
\begin{equation*}
g\in G_{\Phi _{\left( p\right) }^{^{^{\alpha }}}}\Leftrightarrow \left(
E_{n}(g)\right) _{n}\in l_{\Phi _{\left( p\right) }^{^{\alpha
}}}\Leftrightarrow \left( E_{n}(g)\right) _{n}\in l_{\varphi _{\alpha
,p},p}\Leftrightarrow g\in G_{\varphi _{\alpha ,p},p}.
\end{equation*}%
In conclusion 
\begin{equation*}
G_{\Phi _{\left( p\right) }^{^{\alpha }}}=G_{\varphi _{\alpha ,p},p}.
\end{equation*}%
Similarly%
\begin{equation*}
G_{\Phi _{\left( q\right) }^{^{\beta }}}=G_{\varphi _{\beta ,q},q},\text{
respectively }G_{\Phi _{\left( r\right) }^{^{\gamma }}}=G_{\varphi _{\gamma
,r},r}.
\end{equation*}%
The last step is to verify the hypotheses of Theorem\textbf{\ }10\textbf{\ }%
for%
\begin{equation*}
\varphi _{0}:=\varphi _{\alpha ,p},\varphi _{1}:=\varphi _{\beta ,q},\chi
:=f.
\end{equation*}%
From the definition of $f$ we obtain that $\alpha _{\overline{f}}=\beta _{%
\overline{f}}=\frac{1}{l}$ and hence 
\begin{equation*}
0<\alpha _{\overline{f}}=\beta _{\overline{f}}<1.
\end{equation*}%
Because 
\begin{equation*}
\varphi =\frac{\varphi _{0}}{\varphi _{1}}=\frac{\varphi _{\alpha ,p}}{%
\varphi _{\beta ,q}}
\end{equation*}%
we can write 
\begin{equation*}
\overline{\varphi }\left( t\right) =\underset{s>0}{\sup }\,\frac{\varphi
_{\alpha ,p}\left( st\right) }{\varphi _{\beta ,q}\left( st\right) }\cdot 
\frac{\varphi _{\beta ,q}\left( s\right) }{\varphi _{\alpha ,p}\left(
s\right) }.
\end{equation*}%
Now we are interested in computing $\underset{t\rightarrow 0}{\lim }%
\overline{\varphi }\left( t\right) .$ To start we notice that if \ $%
t\rightarrow 0$ then $st<1.$ Consequently,%
\begin{equation*}
\underset{t\rightarrow 0}{\lim }\,\overline{\varphi }\left( t\right) =%
\underset{t\rightarrow 0}{\lim }\,\underset{s>0}{\sup }\,\frac{\varphi
_{\alpha ,p}\left( st\right) }{\varphi _{\beta ,q}\left( st\right) }\cdot 
\frac{\varphi _{\beta ,q}\left( s\right) }{\varphi _{\alpha ,p}\left(
s\right) }=
\end{equation*}%
\begin{equation*}
\underset{t\rightarrow 0}{\lim }\,\underset{s>0}{\sup }\,\frac{\left(
st\right) ^{\frac{1}{p}}}{\left( st\right) ^{\frac{1}{q}}}\cdot \frac{%
\varphi _{\beta ,q}\left( s\right) }{\varphi _{\alpha ,p}\left( s\right) }=%
\underset{t\rightarrow 0}{\lim }\,t^{\frac{1}{p}-\frac{1}{q}}\,\underset{s>0}%
{\sup }\,s^{\frac{1}{p}-\frac{1}{q}}\frac{\varphi _{\beta ,q}\left( s\right) 
}{\varphi _{\alpha ,p}\left( s\right) }.
\end{equation*}%
We only have to evaluate $\underset{s>0}{\sup }\,s^{\frac{1}{p}-\frac{1}{q}}%
\frac{\varphi _{\beta ,q}\left( s\right) }{\varphi _{\alpha ,p}\left(
s\right) }.$ If $s\in (0,1]$ then 
\begin{equation*}
\underset{0<s<1}{\sup }\,s^{\frac{1}{p}-\frac{1}{q}}\frac{\varphi _{\beta
,q}\left( s\right) }{\varphi _{\alpha ,p}\left( s\right) }=\underset{0<s<1}{%
\sup }\,s^{\frac{1}{p}-\frac{1}{q}}\frac{s^{\frac{1}{q}}}{s^{\frac{1}{p}}}=1.
\end{equation*}%
If $s\in (1,\infty )$ then%
\begin{equation*}
\underset{1<s<\infty }{\sup }\,s^{\frac{1}{p}-\frac{1}{q}}\frac{\varphi
_{\beta ,q}\left( s\right) }{\varphi _{\alpha ,p}\left( s\right) }=\underset{%
1<s<\infty }{\sup }\,s^{\frac{1}{p}-\frac{1}{q}}\frac{\left( 1-\left\{
s\right\} \right) \left( \beta _{\left[ s\right] }\left[ s\right] \right) ^{%
\frac{1}{q}}+\left\{ s\right\} \left( \beta _{\left[ s\right] +1}\left( %
\left[ s\right] +1\right) \right) ^{\frac{1}{q}}}{\left( 1-\left\{ s\right\}
\right) \left( \alpha _{\left[ s\right] }\left[ s\right] \right) ^{\frac{1}{p%
}}+\left\{ s\right\} \left( \alpha _{\left[ s\right] +1}\left( \left[ s%
\right] +1\right) \right) ^{\frac{1}{p}}}\leq
\end{equation*}%
\begin{equation*}
\leq \underset{1<s<\infty }{\sup }\,s^{\frac{1}{p}-\frac{1}{q}}\frac{\beta _{%
\left[ s\right] }^{\frac{1}{q}}}{\alpha _{\left[ s\right] +1}^{\frac{1}{p}}}%
\cdot \frac{\left( \left[ s\right] +1\right) ^{\frac{1}{q}}}{\left[ s\right]
^{\frac{1}{p}}}\leq \underset{1<s<\infty }{\sup }\,\left( \frac{\alpha _{%
\left[ s\right] }}{\alpha _{\left[ s\right] +1}}\right) ^{\frac{1}{p}}s^{%
\frac{1}{p}-\frac{1}{q}}\frac{\left( \left[ s\right] +1\right) ^{\frac{1}{q}}%
}{\left[ s\right] ^{\frac{1}{p}}}\leq
\end{equation*}%
\begin{equation*}
\leq M\underset{1<s<\infty }{\sup }\,s^{\frac{1}{p}-\frac{1}{q}}\left[ s%
\right] ^{\frac{1}{q}-\frac{1}{p}}\left( 1+\frac{1}{\left[ s\right] }\right)
^{\frac{1}{q}}\leq M\underset{1<s<\infty }{\sup }\,\left( \left[ s\right]
+1\right) ^{\frac{1}{p}-\frac{1}{q}}\left[ s\right] ^{\frac{1}{q}-\frac{1}{p}%
}\left( 1+\frac{1}{\left[ s\right] }\right) ^{\frac{1}{q}}\leq
\end{equation*}%
\begin{equation*}
\leq M_{1}\underset{1<s<\infty }{\sup }\,\left( 1+\frac{1}{\left[ s\right] }%
\right) ^{\frac{1}{q}}<2^{\frac{1}{q}}M_{1}.
\end{equation*}%
In conclusion for the both cases we obtain $\underset{t\rightarrow 0}{\lim }%
\,\overline{\varphi }\left( t\right) =0$ and furthermore%
\begin{equation*}
\beta _{\overline{\varphi }}>0.
\end{equation*}%
Similarly we analyze the functions $\varphi _{0}$ and $\varphi _{1}.$ We
have 
\begin{equation*}
\overline{\varphi _{\alpha ,p}}\left( t\right) =\underset{s>0}{\sup }\,\frac{%
\varphi _{\alpha ,p}\left( st\right) }{\varphi _{\alpha ,p}\left( s\right) },%
\text{ for every }t\in \left( 0,\infty \right) .
\end{equation*}%
The following inequalitiy is true:%
\begin{equation*}
\overline{\varphi _{\alpha ,p}}\left( t\right) \leq \left( M\left( \frac{1}{t%
}\right) \frac{1}{\left[ \frac{1}{t}\right] }\right) ^{\frac{1}{p}},
\end{equation*}%
hence we can write%
\begin{equation*}
0\leq \underset{t\rightarrow 0}{\lim }\,\overline{\varphi _{\alpha ,p}}%
\left( t\right) \leq \underset{t\rightarrow 0}{\lim }\,\left( M\left( \frac{1%
}{t}\right) \frac{1}{\left[ \frac{1}{t}\right] }\right) ^{\frac{1}{p}}.
\end{equation*}%
But 
\begin{equation*}
\underset{t\rightarrow 0}{\lim }\left( M\left( \frac{1}{t}\right) \frac{1}{%
\left[ \frac{1}{t}\right] }\right) ^{\frac{1}{p}}=0
\end{equation*}%
so we conclude that $\underset{t\rightarrow 0}{\lim }\overline{\varphi
_{\alpha ,p}}\left( t\right) =0$ and furthermore 
\begin{equation*}
\beta _{\overline{\varphi _{0}}}>0.
\end{equation*}%
Similarly we proof that $\underset{t\rightarrow 0}{\lim }\overline{\varphi
_{\beta ,q}}\left( t\right) =0$ from which we derive again%
\begin{equation*}
\beta _{\overline{\varphi _{1}}}>0.
\end{equation*}%
The hypotheses of Theorem 10 being fulfilled we can aply it to obtain%
\begin{equation*}
\left( G_{\varphi _{\alpha ,p},p},G_{\varphi _{\beta ,q},q}\right)
_{f,r}=G_{\rho ,r},
\end{equation*}%
where the function $\rho :\left( 0,\infty \right) \rightarrow \left(
0,\infty \right) $ is defined by%
\begin{equation*}
\rho \left( t\right) :=\frac{\varphi _{\alpha ,p}\left( t\right) }{\left( 
\frac{\varphi _{\alpha ,p}\left( t\right) }{\varphi _{\beta ,q}\left(
t\right) }\right) ^{\frac{1}{l}}}.
\end{equation*}%
We compute%
\begin{equation*}
\rho \left( n\right) =\frac{\varphi _{\alpha ,p}\left( n\right) }{\left( 
\frac{\varphi _{\alpha ,p}\left( n\right) }{\varphi _{\beta ,q}\left(
n\right) }\right) ^{\frac{1}{l}}}=\varphi _{\alpha ,p}\left( n\right) ^{1-%
\frac{1}{l}}\varphi _{\beta ,q}\left( n\right) ^{\frac{1}{l}}=
\end{equation*}%
\begin{equation*}
=\alpha _{n}^{\frac{1}{p}-\frac{1}{p\cdot l}}\beta _{n}^{\frac{1}{ql}}n^{%
\frac{p+ql-q}{pql}}=\left( \alpha _{n}^{r\left( \frac{1}{p}-\frac{1}{pl}%
\right) }\beta _{n}^{\frac{r}{ql}}\right) ^{\frac{1}{r}}n^{\frac{1}{r}%
}=\left( \gamma _{n}n\right) ^{\frac{1}{r}}.
\end{equation*}%
In conclusion $l_{\rho ,r}=l_{\Phi _{\left( r\right) }^{\gamma }},$ where $%
\gamma =\left( \gamma _{n}\right) _{n}$ and hence $G_{\rho ,r}=G_{\varphi
_{\gamma ,r},r}.$ Adding the two equalities proved before $G_{\varphi
_{\alpha ,p},p}=G_{\Phi _{\left( p\right) }^{\alpha }}$ and $G_{\varphi
_{\beta ,q},q}=G_{\Phi _{\left( q\right) }^{\beta }}$ we obtain the desired
result.
\end{proof}

\begin{acknowledgement}
I would like to express my gratitude to Professor N. Tita for many fruitful
discussions and also to Professor C. Merucci for sending me material with
valuable information on "real interpolation method with functional
parameter". Also I want to thank to Professor R. Nest for careful reading of
my manuscript and doing useful remarks and suggestions.
\end{acknowledgement}

\end{document}